\documentclass[12pt]{amsart}
\usepackage{amssymb, amsthm, graphicx, epsfig, graphs, verbatim}
\newtheorem{thm}{Theorem}

\newtheorem{lem}[thm]{Lemma}

\theoremstyle{definition}

\newtheorem{rem}[thm]{Remark}

\newcommand{\ga}{\gamma}

\newcommand{\Si}{\Sigma}

\newcommand{\x}{\times}
\newcommand{\s}{\mathbf s}

\renewcommand{\t}{\mathbf t}

\newcommand{\Z}{\mathbb Z}

\newcommand{\Q}{\mathbb Q}

\newcommand{\co}{\colon\thinspace} 
\newcommand{\hf}{{{\widehat {HF}}}}

\begin{document}
\mathsurround=1pt \title{Heegaard Floer invariants and tight contact
three--manifolds}

\author{Paolo Lisca}
\address{Dipartimento di Matematica\\
Universit\`a di Pisa \\I--56127 Pisa, ITALY} 
\email{lisca@dm.unipi.it}

\author{Andr\'{a}s I. Stipsicz}
\address{R\'enyi Institute of Mathematics\\
Hungarian Academy of Sciences\\
H-1053 Budapest\\ 
Re\'altanoda utca 13--15, Hungary}
\email{stipsicz@math-inst.hu}

\begin{abstract}
Let $Y_r$ be the closed, oriented three--manifold obtained by
performing rational $r$--surgery on the right--handed trefoil knot in
the three--sphere. Using contact surgery and the Heegaard Floer
contact invariants we construct positive, tight contact structures on
$Y_r$ for every $r\neq 1$. This implies, in particular, that the
oriented boundaries of the positive $E_6$ and $E_7$ plumbings carry
positive, tight contact structures, solving a well--known open
problem.
\end{abstract}
\thanks{The first author was partially supported by MURST, and he is a
  member of EDGE, Research Training Network HPRN-CT-2000-00101,
  supported by The European Human Potential Programme. The second
  author would like to thank Peter Ozsv\'ath and Zolt\'an Szab\'o for
  many useful discussions regarding their joint work. The second
  author was partially supported by OTKA T034885}

\maketitle

\section{Introduction}\label{s:intro}

One of the central problems in contact topology is the construction of
positive, tight contact structures on closed, oriented
three--manifolds. Using contact surgery~\cite{DG1, DG2} one can easily
find contact structures on various three--manifolds, but until
recently it has been a hard problem to check that such structures are
tight. The purpose of this paper is to show that the newly introduced
Heegaard Floer contact invariants~\cite{OSz6} provide an appropriate
tool to attack this problem. In fact, under favourable circumstances,
a partial determination of the Heegaard Floer contact invariants
suffices to prove that the structure under examination is tight. Here
we present a prototype of the results obtainable along these lines. A
more thorough investigation will appear later~\cite{LS}.

\begin{thm} \label{t:main}
Let $r\in\Q\cup\{\infty\}$, and denote by $Y_r$ the closed, oriented
three--manifold obtained by performing $r$-surgery on the right-handed
trefoil knot $K\subset S^3$. If $r\neq 1$, then $Y_r$ carries a
positive, tight contact structure.
\end{thm}

{\bf Remarks.  (a)} Theorem~\ref{t:main} solves, in particular, the
well--known open problem asking whether the boundaries of the positive
definite $E_6$ and $E_7$ plumbings (respectively $Y_3$ and $Y_2$ in
the above notation) carry positive, tight contact
structures~\cite{EN}. Note that the contact structures on $Y_2$ and
$Y_3$ given by Theorem~\ref{t:main} are not symplectically
semi--fillable because the three--manifolds $Y_2$ and $Y_3$ do not
carry symplectically semi--fillable contact structures~\cite{Paolo}.

{\bf (b)} The three--manifold $Y_1$ is the oriented boundary of the
positive definite $E_8$ plumbing. It is known that $Y_1$ does not
carry positive, tight contact structures~\cite{EH}.

In proving Theorem~\ref{t:main} we shall use contact surgery to define
contact structures on $Y_r$ ($r\neq 1$) and then show that the
Heegaard Floer contact invariants of those structures do not vanish,
implying tightness. During the course of the proof we show that the
contact invariants are nontrivial for infinitely many tight, not
semi--fillable contact three--manifolds. 

\section{The proof of Theorem~\ref{t:main}}\label{s:proof}

\begin{figure}[ht]
\begin{center}
\epsfig{file=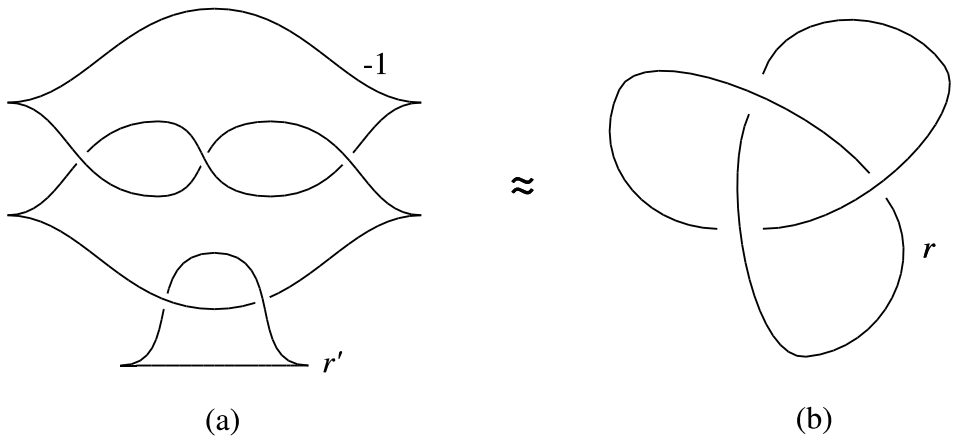, height=6cm}
\end{center}
\caption{Contact structures on $Y_r$ with $r=\frac{1}{1-r'}$}
\label{f:structure}
\end{figure}

Consider the contact structures defined by Figure~\ref{f:structure}(a).
The meaning of such a diagram is the following: perform $(-1)$-- and
$r'$--surgeries on the Legendrian knots of the diagram with respect to
the contact framings. Extend the standard contact structure of $S^3$
to the new solid tori by structures which are tight on these
tori. Such an extension exists if and only if $r'\neq 0$ and is unique
if and only if $r'=\frac{1}{k}$ for some integer $k\in \Z$. In
general, the number of different extensions depends on a continued
fraction expansion of $r'$~\cite{DG1, DG2}. As illustrated in
Figure~\ref{f:structure}(b), this procedure defines contact structures
on $Y_r$ for $r=\frac{1}{1-r'}$ with $r'\neq 0$, i.e., for all $r\neq
1$.

\begin{figure}[ht]
\begin{center}
\epsfig{file=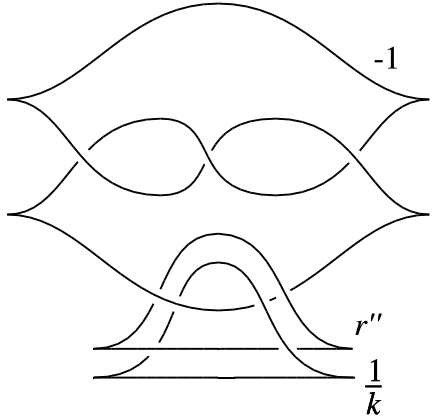, height=6cm}
\end{center}
\caption{Equivalent surgery diagram for contact structures on $Y_r$}
\label{f:modify}
\end{figure}

According to~\cite{DG2}, Proposition~3, a contact $r'$--surgery with
$r'<0$ or $r'=\infty$ can be replaced by a sequence of Legendrian
(i.e., contact $(-1)$--) surgeries, providing a Stein fillable and
hence tight contact three--manifold. If $r'>0$, $r'\neq\infty$, then
by~\cite{DG2}, Proposition~7, the contact structures of
Figure~\ref{f:structure}(a) can be equivalently given by the diagram
of Figure~\ref{f:modify} for any positive integer $k$ and
$r''=\frac{r'}{1-kr'}$.

In a series of papers~\cite{OSzF1, OSzF2, OSzcobor, OSz6} Ozsv\'ath and
Szab\'o introduced a family of invariants associating a homology group
$\hf(M,\s)$ to a Spin$^c$ three--manifold $(M,\s)$, a homomorphism
\[
F_{W,\t}\colon \hf (M_1, \s _1) \to \hf (M_2, \s _2)
\] 
to a Spin$^c$ cobordism $(W,\t)$ and an element 
\[
c(M, \xi )\in \hf (-M, \s _{\xi})/\langle \pm 1 \rangle
\] 
to a contact three--manifold $(M,\xi )$, where $\s _{\xi }$ denotes
the Spin$^c$ structure induced by the contact structure $\xi$. In this
note we shall always use this homology theory with $\Z/2\Z$
coefficients, so that the above sign ambiguity for $c(M,\xi)$ does not
occur.  It is proved in~\cite{OSz6} that if $(M, \xi)$ is overtwisted
then $c(M,\xi)=0$, and if $(M,\xi)$ is Stein fillable then
$c(M,\xi)\neq 0$. The following lemma is implicitely contained in
Theorem~4.2 of~\cite{OSz6}.

\begin{lem}\label{l:property3}
Suppose that $(M_2,\xi_2)$ is obtained from $(M_1,\xi _1)$ by a
contact $(+1)$--surgery. Then we have 
\[
F_{-X} (c(M_1, \xi_1))= c(M_2,\xi_2), 
\]
where $-X$ is the cobordism induced by the surgery with reversed
orientation and $F_{-X}$ is the sum of $F_{-X,\s}$ over all Spin$^c$
structures $\s$ extending the Spin$^c$ structures on $-M_i$ by
$\xi_i$, $i=1,2$.  In particular, if $c(M_2, \xi_2)\neq 0$ then $(M_1,
\xi_1)$ is tight.
\end{lem}

\begin{proof}
Let us assume that we are performing contact $(+1)$--surgery along the
Legendrian knot $\ga\subset (M_1,\xi_1)$. It follows from~\cite{AO}
and~\cite{DG2} that there is an open book decomposition $(F,\phi)$ on
$M_1$ compatible with $\xi_1$ in the sense of Giroux and such that
$\ga$ lies on a page and is not homotopic to the boundary. Then, an
open book for $(M_2,\xi_2)$ is given by $(F,\phi')$, where
$\phi'=\phi\circ R^{-1}_\ga$ and $R_\ga$ is the right--handed Dehn
twist along $\ga$. The first part of the statement now follows
applying Theorem~4.2 of~\cite{OSz6}. The second part of the statement
follows immediately from the fact that the invariant of an overtwisted
contact structure vanishes.
\end{proof}

We will now describe an extremely effective computational tool for
Heegaard Floer homology, the~\emph{surgery exact triangle}. Let $M$ be
a closed, oriented three--manifold and let $K\subset M$ be a framed
knot. For $n \in\Z$, let $M_n$ denote the three--manifold given by
$n$--surgery along $K\subset M$ with respect to the given framing, and
call the resulting cobordism $X_n$. For a given Spin$^c$ structure
$\s$ on $M\setminus K$ let
\[
\hf (M_n , [\s ]) = \bigoplus _{\t \in {\mathfrak {S}}}\hf (M_n , \t), 
\]
where ${\mathfrak {S}}=\{ \t \in Spin ^c (M_n)\mid \t \vert
_{M_n\setminus K}= \s \}$, and define
\[
\hf(M_n) = \bigoplus_{\t\in Spin ^c (M_n)} \hf (M_n , \t).
\] 
The cobordism $X_n$ induces a homomorphism
\[
F_{X_n}\co \hf (M, [\s] )\to \hf (M_n , [\s ])
\]
obtained by summing over Spin$^c$ structures. 

\begin{thm} [\cite{OSzF2}, Theorem~9.16]\label{t:triangle}
The homomorphism $F_{X_n}$ fits into an exact triangle:
\[
\begin{graph}(6,2)
\graphlinecolour{1}\grapharrowtype{2}
\textnode {A}(1,1.5){$\hf (M, [\s] )$}
\textnode {B}(5, 1.5){$\hf (M_n , [\s ])$}
\textnode {C}(3, 0){$\hf (M_{n+1}, [\s ])$}
\diredge {A}{B}[\graphlinecolour{0}]
\diredge {B}{C}[\graphlinecolour{0}]
\diredge {C}{A}[\graphlinecolour{0}]
\freetext (3,1.8){$F_{X_n}$}
\end{graph}
\]
\qed\end{thm}

\noindent
{\bf Remark}. Since the Spin$^c$ structures on $M_n$ are partitioned 
according to their restrictions to $M_n\setminus K$, by taking the direct
sum of the exact triangles given by Theorem~\ref{t:triangle} one obtains
an analogous exact triangle involving the groups $\hf(M)$, $\hf(M_n)$ and 
$\hf(M_{n+1})$. 

\begin{lem}\label{l:k=1}
Let $K\subset (S^3,\xi_{st})$ be the Legendrian unknot with
Thurston--Bennequin invariant equal to $-1$ and vanishing rotation
number. Then, the contact three--manifold $(S^1\times S^2, \eta)$
obtained by $(+1)$--surgery along $K$ has nonvanishing Heegaard Floer
invariant.
\end{lem}

\begin{proof}
In view of Lemma~\ref{l:property3}, Theorem~\ref{t:triangle} and 
the remark following the latter, we have an exact triangle: 
\[
\begin{graph}(6,2)
\graphlinecolour{1}\grapharrowtype{2}
\textnode {A}(1,1.5){$\hf(S^3)$}
\textnode {B}(5, 1.5){$\hf(S^1\x S^2)$}
\textnode {C}(3, 0){$\hf(S^3)$}
\diredge {A}{B}[\graphlinecolour{0}]
\diredge {B}{C}[\graphlinecolour{0}]
\diredge {C}{A}[\graphlinecolour{0}]
\freetext (2.8,1.8){$F_{-X}$}
\end{graph}
\]

where $-X$ is the cobordism from $S^3$ to $S^1\x S^2$ obtained 
by attaching a two--handle to $S^3$ along a zero--framed unknot, 
and 
\[
F_{-X}(c(S^3, \xi_{\rm st}))=c(S^2\x S^1, \eta).
\]
By~\cite{OSzF2}, $\hf(S^1\x S^2)$ is isomorphic to $(\Z/2\Z)^2$, while
$\hf(S^3)$ is isomorphic to $\Z/2\Z$. Exactness of the triangle
immediately implies that $F_{-X}$ is injective. Since $(S^3,\xi_{st})$
is Stein fillable we have $c(S^3,\xi_{st})\neq 0$, therefore $c(S^2\x
S^1, \eta)\neq 0$.
\end{proof}

\begin{figure}[ht]
\begin{center}
\epsfig{file=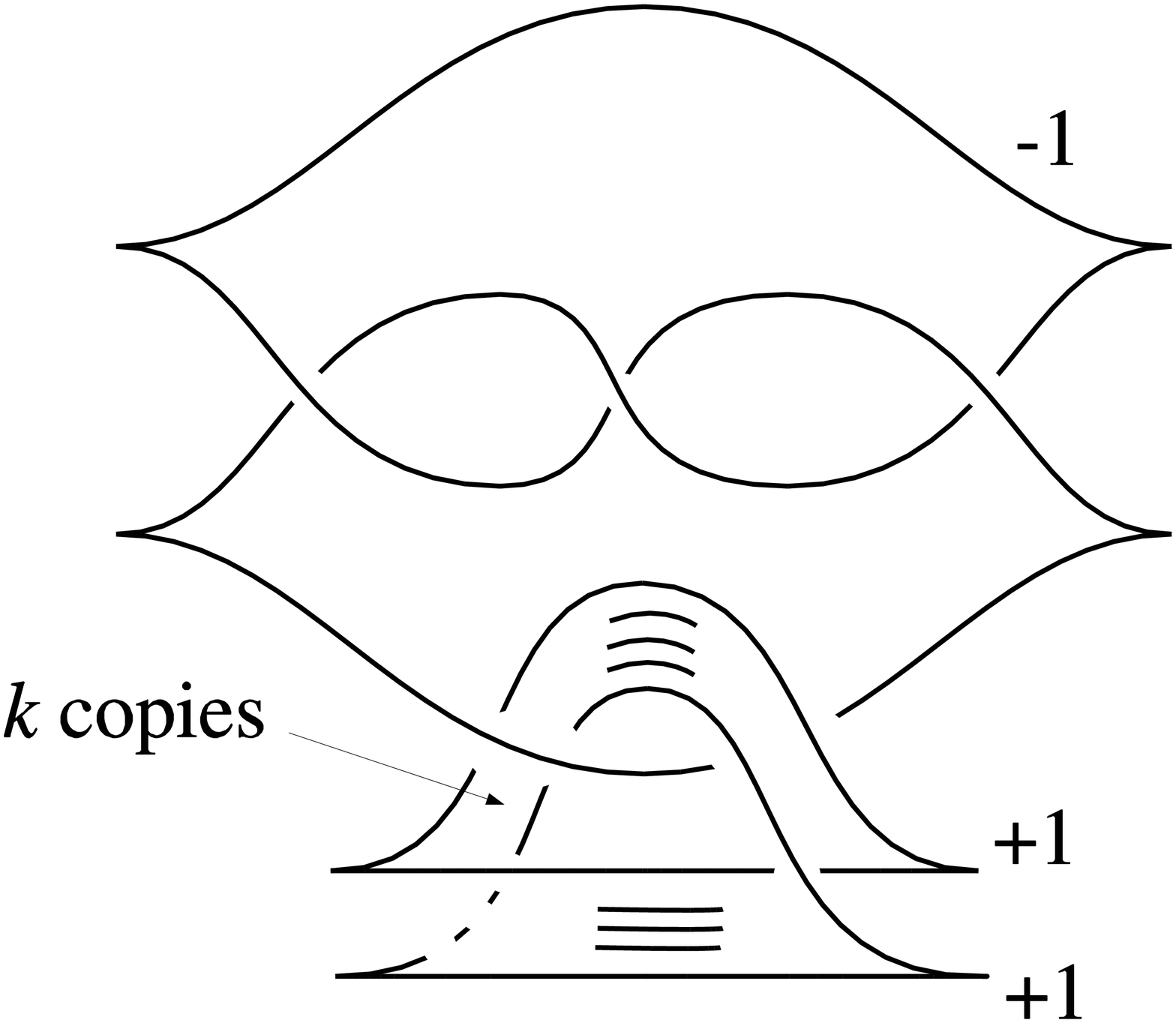, height=6cm}
\end{center}
\caption{Equivalent surgery diagram for $(V_k,\xi_k)$}
\label{f:equiv}
\end{figure}

Let $(V_k, \xi _k)$ denote the contact three--manifold obtained by
choosing $r'=\frac{1}{k}$ in Figure~\ref{f:structure}(a), so that $V_k
\cong Y_{\frac{k}{k-1}}$.  Notice that for $r'=k=1$ the
three--manifold $V_1\cong Y_{\infty} $ is diffeomorphic to the
three--sphere $S^3$. According to~\cite{DG2} the contact
three--manifold $(V_k, \xi _k)$ can be alternatively defined by the
diagram of Figure~\ref{f:equiv}, which contains $k$ contact
$(+1)$--framed Legendrian unknots.

\begin{lem}\label{l:Vk}
Let $k\geq 1$ be an integer. Then, $\dim_{\Z/2\Z}\hf(-V_k)=k$.
\end{lem}

\begin{proof}
$(V_{k+1}, \xi_{k+1})$ is clearly obtained from $(V_k,\xi_k)$ by
performing a contact $(+1)$--surgery. Therefore, as in the proof of
Lemma~\ref{l:k=1} the cobordism $X_k$ corresponding to the surgery
induces a homomorphism $F_{-X_k}$ which fits into an
exact triangle:
\[
\begin{graph}(8,2)
\graphlinecolour{1}\grapharrowtype{2}
\textnode {A}(3,1.5){$\hf (-V_k)$}
\textnode {B}(7, 1.5){$\hf (-V_{k+1})$}
\textnode {C}(5, 0){$\hf (-Y_{+1})$}
\diredge {A}{B}[\graphlinecolour{0}]
\diredge {B}{C}[\graphlinecolour{0}]
\diredge {C}{A}[\graphlinecolour{0}]
\freetext (4.8,1.8){$F_{-X_k}$}
\freetext(0, 1){$(*)$}
\end{graph}
\]

A simple computation shows that $Y_{+1}$ is the Poincar\'e
sphere $\Si(2,3,5)$, and it follows from the calculations
of~\cite{OSzabs}, Section~3.2, that $\hf (-Y_{+1})=\Z/2\Z$. Therefore,
setting $d(k)=\dim_{\Z/2\Z}\hf(-V_k)$, Triangle $(*)$ implies
\begin{equation}\label{e:equal}
d(k+1) = d(k) \pm 1
\end{equation}
for every $k\geq 1$. Now observe that $-V_k$ can be presented by the
surgery diagram of Figure~\ref{f:plum}. 
\begin{figure}[ht]
\begin{center}
\epsfig{file=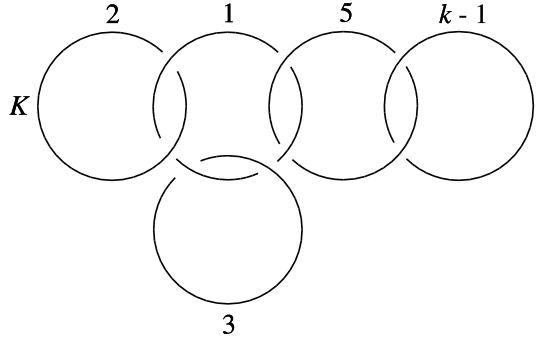, height=5cm}
\end{center}
\caption{A surgery diagram for $-V_k$}
\label{f:plum}
\end{figure}
Consider the three--manifold $M$ obtained by surgery on the framed
link of Figure~\ref{f:plum} with the $2$--framed knot $K$
deleted. Applying Theorem~\ref{t:triangle} to the pair $(M,K)$ for
$n=1$ yields the exact triangle:
\[
\begin{graph}(6,2)
\graphlinecolour{1}\grapharrowtype{2}
\textnode {A}(1,1.5){$\hf(L(7k-9,7))$}
\textnode {B}(5, 1.5){$\hf (L(8k-9,8))$}
\textnode {C}(3, 0){$\hf (-V_k)$}
\diredge {A}{B}[\graphlinecolour{0}]
\diredge {B}{C}[\graphlinecolour{0}]
\diredge {C}{A}[\graphlinecolour{0}]
\end{graph}
\]

Since by~\cite{OSzF2}, Proposition~3.1, $\dim_{\Z/2\Z}\hf(L(p,q))=p$
for every $p$ and $q$, exactness of the triangle implies
\begin{equation}\label{e:ineq}
d(k) \geq k
\end{equation}
for every $k\geq 1$, and since $V_1$ is diffeomorphic to $S^3$, we have
$d(1)=1$. Therefore, by~\eqref{e:equal} and~\eqref{e:ineq} we have
$d(k)=k$ for every $k\geq 1$.
\end{proof}

\begin{lem}\label{l:nonvanish}
Let $k\geq 1$ be an integer. Then, $c(V_k, \xi_k)\neq 0$.
\end{lem}
\begin{proof}
Consider Figure~\ref{f:equiv} for $k=1$, which represents
$(V_1,\xi_1)$. By~\cite{DG1}, performing a contact $(+1)$--surgery on
a Legendrian pushoff of the Legendrian trefoil is equivalent to
erasing the trefoil from the picture, thus resulting in $(S^1\times
S^2, \eta)$. It follows from Lemmas~\ref{l:property3} and~\ref{l:k=1}
that $c(V_1,\xi _1)\neq 0$. Since by~\cite{OSzabs} we know that $\hf
(-Y_{+1})=\Z/2\Z$, it follows from Lemma~\ref{l:Vk} and Triangle~$(*)$
that the homomorphism $F_{-X_k}$ is injective for every
$k\geq 1$. Thus,
\[
c(V_k, \xi _k)=F_{-X_{k-1}}(c(V_{k-1},\xi_{k-1}))\neq 0
\]
for every $k\geq 1$.
\end{proof}

{\bf Remark}. As we have already remarked, $(V_2,\xi_2)$ is a tight
but not semi--fillable contact three--manifold. Since contact
$(+1)$--surgery on a nonfillable structure produces a nonfillable
structure~\cite{DG1, EH1}, by Lemma~\ref{l:nonvanish} the contact
three--manifold $(V_k, \xi_k)$ is tight, not symplectically
semi--fillable for any $k\geq 2$.

\begin{proof}[Proof of Theorem~\ref{t:main}]
If $r'<0$ or $r'=\infty$, any contact surgery given by
Figure~\ref{f:structure}(a) can be realized by a Legendrian surgery,
therefore the resulting contact structure is Stein fillable and hence
tight. If $r'\neq\infty$ and $r'>0$, choose an integer $k$ so large
that $r''=\frac{r'}{1-kr'} < 0$. Since contact $(+1)$--surgery on a
Legendrian pushoff cancels Legendrian surgery, Figure~\ref{f:modify}
shows that one can perform a sequence of contact $(+1)$--surgeries on
any contact structure given by the diagram of
Figure~\ref{f:structure}(a), obtaining $(V_k, \xi_k)$. It follows from
Lemma~\ref{l:nonvanish} and repeated applications of
Lemma~\ref{l:property3} that all the contact structures defined by
Figure~\ref{f:structure}(a) are tight.
\end{proof}


\begin{thebibliography}{AAA}

\bibitem{AO}
{\bf S.~Akbulut and  B.~Ozbagci},
{\it Lefschetz fibrations on compact Stein surfaces},   
Geom. Topol.  {\bf 5}  (2001), 319--334.

\bibitem{DG1}
{\bf F.~Ding and H.~Geiges},
{\it Symplectic fillability of tight contact structures on 
torus bundles}, Alg. and Geom. Topol. {\bf1} (2001), 153--172.

\bibitem{DG2} {\bf F.~Ding and H.~Geiges}, 
{\it A Legendrian surgery presentation of contact 3-manifolds}, 
arXiv preprint math.SG/0107045, to appear on 
Math. Proc. Cambridge Philos. Soc.

\bibitem{EH}
{\bf J.~Etnyre and K.~Honda},
{\it On the nonexistence of tight contact structures}
Ann. of Math. {\bf153} (2001), 749--766.

\bibitem{EH1}
{\bf J.~Etnyre and K.~Honda},
{\it Tight contact structures with no symplectic fillings},
Invent. Math. {\bf 148} (2002), no. 3, 609--626.

\bibitem{EN}
{\bf J.~Etnyre and L.~Ng},
{\it Problems in Low Dimensional Contact Topology}, 
to appear in Proceedings of the 2001 Georgia International 
Geometry and Topology Conference.

\bibitem{PLpos}
{\bf P.~Lisca},
{\it Symplectic fillings and positive scalar curvature},
Geom. Topol. {\bf2} (1998) 103--116.

\bibitem{Paolo}
{\bf P.~Lisca},
{\it On symplectic fillings of 3-manifolds},
Proceedings of the $6^{th}$ G\"okova Geometry-Topology
Conference, Turkish J. Math. {\bf23} (1999) 151--159.

\bibitem{LS}
{\bf P.~Lisca and A.~Stipsicz},
{\it  Surgered contact structures and Heegaard Floer invariants},
in preparation.

\bibitem{OSzF1}
{\bf P.~Ozsv\'ath and Z.~Szab\'o},
{\it Holomorphic disks and topological invariants for
rational homology three-spheres},
to appear in Ann. of Math. (arXiv:math.SG/0101206)

\bibitem{OSzF2}
{\bf P.~Ozsv\'ath and Z.~Szab\'o},
{\it Holomorphic disks and three-manifold invariants: 
properties and applications},
to appear in Ann. of Math. (arXiv:math.SG/0105202)

\bibitem{OSzcobor}
{\bf P.~Ozsv\'ath and Z.~Szab\'o},
{\it Holomorphic triangles and invariants for smooth four-manifolds}, 
arXiv preprint math.SG/0110169. 

\bibitem{OSzabs}
{\bf P.~Ozsv\'ath and Z.~Szab\'o},
{\it Absolutely Graded Floer homologies and intersection forms for
four-manifolds with boundary}, arXiv preprint math.SG/0110170. 

\bibitem{OSz6}
{\bf P.~Ozsv\'ath and Z.~Szab\'o},
{\it Heegaard Floer homologies and contact structures},
arXiv preprint math.SG/0210127.


\end{thebibliography}
\end{document}